%

\input amstex

\magnification 1200
\loadmsbm
\parindent 0 cm

\define\nl{\bigskip\item{}}
\define\snl{\smallskip\item{}}
\define\inspr #1{\parindent=20pt\bigskip\bf\item{#1}}
\define\iinspr #1{\parindent=27pt\bigskip\bf\item{#1}}
\define\einspr{\parindent=0cm\bigskip}

\define\co{\Delta}

\define\ot{\otimes}

\define\tr{\triangleright}
\define\tl{\triangleleft}

\input amssym
\input amssym.def

\centerline{\bf Tools for working with multiplier Hopf algebras}
\bigskip\bigskip
\centerline{\it  A.\ Van Daele \rm ($^*$)}
\bigskip\bigskip\bigskip
{\bf Abstract} 
\bigskip 
Let $(A,\co)$ be a multiplier Hopf algebra. In general, the underlying algebra $A$ need not have an identity and the coproduct $\Delta$ does not map $A$ into $A\ot A$ but rather into its multiplier algebra $M(A\ot A)$. In this paper, we study {\it some tools} that are frequently used when dealing with such multiplier Hopf algebras and that are typical for working with algebras without identity in this context.
\snl
The {\it basic ingredient} is  a unital left $A$-module $X$. And the basic construction is that of extending the module by looking at linear maps $\rho:A\to X$ satisfying $\rho(aa')=a\rho(a')$ where $a,a'\in A$. We write the module action as multiplication. Of course, when $x\in X$, and when $\rho(a)=ax$, we get such a linear map. And if $A$ has an identity, all linear maps $\rho$ have this form for $x=\rho(1)$. However, the point is that in the case of a non-unital algebra, the space of such maps is in general strictly bigger than $X$ itself. We get an {\it extended module}, denoted by $\overline X$ (for reasons that will be explained in the paper).
\snl
We study all sorts of more complicated situations where such extended modules occur and we illustrate all of this with {\it several examples}, from very simple ones to more complex ones where iterated extensions come into play. We refer to cases that appear in the literature.
\snl
We use this basic idea of extending modules to explain, in a more rigorous way,  the so-called {\it covering technique}, which is needed when using {\it Sweedler notations} for coproducts and coactions. Again, we give many examples and refer to the existing literature where this technique is applied. 
\nl
\nl
{\it June 2008} (Version 1.1)
\vskip 1 cm
\hrule
\bigskip\parindent 0.6 cm
\item{($^*$)} Department of Mathematics, K.U.\ Leuven, Celestijnenlaan 200B,
B-3001 Heverlee (Belgium). E-mail: Alfons.VanDaele\@wis.kuleuven.be
\parindent 0 cm
 
\newpage

\bf 0. Introduction \rm
\nl
If $A$ is a {\it finite-dimensional Hopf algebra}, the dual space $A'$ of all linear functionals on $A$ is again a Hopf algebra. The product and the coproduct on $A'$ are obtained by dualizing the coproduct and the product of $A$. If $A$ is an {\it infinite-dimensional} Hopf algebra, the result is no longer valid. In some cases, it is possible to find one or more Hopf algebra's $B$ that can be paired with $A$, in a non-degenerate way, and such that the product and coproduct on $B$ are dual to the coproduct and product on $A$.
\snl
On the other hand, if $A$ is any {\it Hopf algebra with integrals} (and invertible antipode), there is a natural subspace $\widehat A$ of the full dual space $A'$, still separating points of $A$, and carrying the structure of a {\it multiplier Hopf algebra}. The product on $\widehat A$ is inherited from the product on $A'$, which is dual to the coproduct on $A$. Also the coproduct $\widehat\Delta$ on $\widehat A$ is dual to the product on $A$ but it does not map $\widehat A$ into the tensor product $\widehat A\ot \widehat A$, but rather into its multiplier algebra $M(\widehat A\ot \widehat A)$. Also this dual multiplier Hopf algebra $(\widehat A,\widehat\Delta)$ has integrals.
\snl
This duality is a special case of the more general duality of regular multiplier Hopf algebras with integrals, the so-called {\it algebraic quantum groups} (see [VD2]. It extends the duality of finite-dimensional Hopf algebras and many of the results about Hopf algebras that can only be proven in the finite-dimensional case, are still true for this much bigger class of multiplier Hopf algebras with integrals. This duality also has some advantages over that of a pairing of Hopf algebras because of the presence of integrals on both sides.
\snl
The {\it main technical difficulty} is due to the fact that the algebras involved, in general, have no identity. In this context, it is also not a good idea to simply add an identity. On the contrary, the relevant construction is that of passing to the {\it multiplier algebra} $M(A)$ when $A$ had no identity. This multiplier Hopf algebra can be characterized as the largest algebra with identity, containing $A$ as a dense two-sided ideal. It is called a dense ideal because of the property that, when $x\in M(A)$, then $x=0$ if either $ax=0$ for all $a\in A$ or $xa=0$ for all $a\in A$.
\snl
This fact has many consequences and implies that often some care is necessary. A good example to illustrate this statement is the {\it use of the Sweedler notation} for multiplier Hopf algebras. It is well established for Hopf algebras and very useful because it makes formula's and calculations more transparent and therefore easier to understand. This is also true when working with multiplier Hopf algebras, but then the use of the Sweedler notation is much more tricky. The (obvious) reason is that the coproduct does not have range in the tensor product, but rather in the multiplier algebra of the tensor product.
\snl
The use of the Sweedler notation for multiplier Hopf algebras has been introduced first in [Dr-VD]. Since then it has been a common practice. One always has to make sure that the {\it proper coverings} are present (again see [Dr-VD]). Recently however, in an attempt to generalize Majid's bicrossproduct construction, we felt the need to develop this tool in a {\it more rigorous way}. Indeed, formula's in the context of these bicrossproducts, where there is a mixing of coproducts and coactions, are quite complicated. The use of the Sweedler notation is a clarifying tool, but quite involved when generalizing to multiplier Hopf algebras. 
\nl
This takes us to the {\it aim} and the {\it content of this paper}.
\nl
In {\it Section} 1 we begin with our {\it basic construction}. We consider an algebra $A$ with a non-degenerate product and a non-degenerate left $A$-module $X$. The action is written as multiplication. We define the space $Y$ of linear maps $\rho:A\to X$ satisfying $\rho(aa')=a\rho(a')$ for all $a,a'\in A$. The space $Y$ is again a left $A$-module in a natural way and it contains $X$ as a submodule. For all elements $y$ in $Y$ we have $ay\in X$ for all $a\in A$. Roughly speaking, this property characterizes $Y$. We call $Y$ the {\it completed module}. Observe that we get nothing new when $A$ has an identity, but if this is not the case, we do get a strictly bigger module in general.
\snl
Also in Section 1, we consider {\it right modules}, {\it bimodules} and we apply the basic construction to {\it more complicated} and {\it iterated cases}. We illustrate all of this with many examples. Most of these come from the theory of multiplier Hopf algebras.
\snl
In {\it Section} 2, we start with the ingredients of the previous section and we consider linear maps $F$ from the left $A$-module $X$ to any vector space $V$. We say that the variable is {\it covered} if there exists an element $e$ in $A$ such that $F(ex)=F(x)$ for all $x$ in $X$. Then, we can {\it extend} the map $F$ to the completed module $Y$ such that $F(y)=F(ey)$ with $e$ as before. Again we consider the case of a right module and a bimodule as well and we give plenty of examples, some of which involving the more complicated iterated cases as given in the first section.
\snl
Also in this section, we discuss the use of the {\it Sweedler notation}. As it turns out, a rigorous treatment is based on the previously introduced notions of covering and the extension of linear maps from the module to the extended module. Also here we give several examples.
\snl
It can be noticed already in the two previous sections that the theory is rather simple and not very deep but that nevertheless, the applications and the examples can become quite complicated. There are also many different situations where the theory can be applied. For this reason, we have included a {\it third section} with more complicated cases and more  examples.
\snl
The {\it last} (very short) {\it section} contains some conclusions.
\snl
Finally, in a (also short) {\it appendix}, we discuss a topological point of view. It is not really used in our approach, but we feel it is illuminating and can help for a better and more complete understanding of the matter.
\snl
The {\it aim of the paper} mainly is to clarify the use of the Sweedler notation, both for coproducts and coactions, by making it more rigorous. The focus in the paper is on the examples as explained before. And because the Sweedler notation has already been used before in this context, we will often give examples that refer to the existing cases in literature. This may give the (wrong) impression that the arguments in earlier work are not completely correct. Indeed, this is not the case. We just explain here, in a more explicit way, the ideas that are behind the techniques used in earlier work. We have noticed that for many people reading this work, the passage from Hopf algebras to multiplier Hopf algebras is not an obvious one. With this paper, we hope to contribute to making this step easier.
\nl
\it Notations and conventions \rm
\nl
We work with (associative) algebras over the field $\Bbb C$ of complex numbers but in fact, the theory is also valid for algebras over other fields. We do not require the algebras to have an identity, but we do assume that always the {\it product is non-degenerate} as a bilinear map. In that case, one can consider the multiplier algebra $M(A)$ as we  explained already before in this introduction. If already $A$ has an identity, then $M(A)=A$, but otherwise, in general, $M(A)$ will be much larger than $A$ itself.  
\snl
We use $1$ for the identity in $M(A)$. On the other hand, we always use $\iota$ to denote the identity map. When e.g.\ we have a linear map $f:A\to \Bbb C$, we use $f\ot\iota$ for the {\it slice map} from $A\ot A$ to $A$ sending $a\ot a'$ to $f(a)a'$.
\snl
As it should be clear from this introduction, we will also use the Sweedler notation, both for coproducts and for actions. But then, we will explain carefully how this should be done.
\snl
The new notations we will introduce further in this paper are chosen very carefully. They may seem somewhat confusing in the first place, but soon, the reader will find out that they are quite logical and convenient. This practice will be completely in agreement with the aim and the spirit of the paper as explained before.
\nl
\it Basic references \rm
\nl
For the basic theory of Hopf algebras, we refer to the standard works [A] and [S]. For multiplier Hopf algebras, the main reference is [VD1] and when they have integrals, it is [VD2]. Actions of multiplier Hopf algebras have been studied in [Dr-VD-Z].  The use of the Sweedler notation for multiplier Hopf algebras has been introduced in [Dr-VD]. See also [VD-Z] for a survey on the subject. The paper [De-V-W], where Majid's bicrossproduct construction is generalized to multiplier Hopf algebras, has been the motivating source for writing this note.
\nl\nl
\bf Acknowledgements \rm
\nl
I would like to thank my coworkers in this field. In particular, I have appreciated very much the hospitality of Shuanhong Wang (and his research group) during my recent visit to China (March 2008).
\nl\nl

\bf 1.  Module extensions \rm
\nl
Let $A$ be an (associative) algebra over the field $\Bbb C$ of complex numbers. It may or may not have an identity. If the algebra has an identity, we denote it by $1$. If the algebra as no identity, we want the product to be non-degenerate (as a bilinear form). Remark that this is automatic if the algebra has an identity.
\snl
Let $X$ be a vector space over $\Bbb C$. Assume that it is a {\it left $A$-module}. We think of the module action as multiplication and therefore we write $ax$ where $a$ is in $A$ and acts on the vector $x\in X$. So, the module property writes as associativity: For all $a,a'\in A$ and $x\in X$, we have $(aa')x=a(a'x)$.
\snl
We will also work with right modules and bimodules (with left and right action, either of the same algebra or two different algebras). And we will use similar notations as in the case of left modules.
\snl
We will always assume {\it non-degenerate module actions}. So, if e.g.\ $X$ is a left $A$-module, we require that $x=0$ if $x\in X$ and $ax=0$ for all $a\in A$. Sometimes, we will assume that the actions are unital. Again for a left $A$-module $X$ this means that the span of the elements $ax$ with $a\in A$ and $x\in X$ is all of $X$. When $A$ has an identity, the action is non-degenerate if and only if $1x=x$ for all $x\in X$. Then it is also unital. When $A$ has no identity, the two notions may be different. For most of the algebras we will consider here, any unital action is also non-degenerate (but also then, the converse need not be true). 
\nl
\it The basic construction \rm
\nl
Consider an algebra $A$ and a left $A$-module $X$. So we assume that the product in $A$ is non-degenerate and that also the module action is non-degenerate. We will now extend the module $X$. The construction is quite simple, but it is the basis of more complicated constructions we will have in the sequel.

\inspr{1.1} Definition \rm 
We denote by $Y$ the vector space of linear maps $\rho:A\to X$ satisfying $\rho(aa')=a\rho(a')$ for all $a,a'\in A$. 
\einspr

For any $x\in X$, we have an element $\rho_x$ in $Y$ defined by $\rho_x(a)=ax$. The map $x\mapsto\rho_x$ is injective. Indeed, if $\rho_x=0$, then $ax=0$ for all $a\in A$ and as the module is assumed to be non-degenerate, we get $x=0$.
\nl
We can make $Y$ into a left $A$-module and we get the following result:

\inspr{1.2} Proposition \rm 
For any $\rho\in Y$ and $a\in A$, define $a\rho:A\to X$ by $(a\rho)(a')=\rho(a'a)$. This makes $Y$ into a left $A$-module, containing $X$ as a sub-module.

\snl\bf Proof\rm:
First, take $a,a',a''\in A$. Then
$$(a\rho)(a'a'')=\rho(a'a''a)=a'\rho(a''a)=a'(a\rho)(a'').$$
So, $a\rho\in Y$. Moreover
$$\align  ((aa')\rho)(a'')&=\rho(a''(aa'))=\rho((a''a)a')\\
                          &=(a'\rho)(a''a)=(a(a'\rho))(a'')
\endalign$$  
so that $(aa')\rho=a(a'\rho)$. It follows that $Y$ is a  left $A$-module. Finally, to show that $X$ is a sub-module, take also $x\in X$ and consider the map $x\mapsto \rho_x$ as before. Then 
$$\rho_{ax}(a')=a'ax=\rho_x(a'a)=(a\rho_x)(a')$$
and we see that $x\mapsto \rho_x$ is a map of left modules. \hfill $\square$
\einspr

If the algebra $A$ has an identity, then  $Y=X$ because when $\rho\in Y$ and $x=\rho(1)$, we get 
$$\rho(a)=a\rho(1)=ax=\rho_x(a).$$ 
In general however, we will have that $X$ is a proper subset of $Y$.
\snl
Also in general, the extended module can be degenerate. Indeed, suppose that $\rho\in Y$ and $a\rho=0$ for all $a\in A$. Then $\rho(a'a)=0$ for all $a,a'\in A$. In order to conclude that $\rho=0$, it is sufficient to have that $A^2=A$. Then, the extended module will be non-degenerate. Under this assumption, when we repeat the procedure, we will end up with the same extended module $Y$.
\snl
We always have that $AY\subseteq X$. If the left module $X$ is {\it unital}, i.e.\ when $X=AX$, then $X\subseteq AY\subseteq X$ so that $AY=X$. In general however, we may have that $AY$ is a proper submodule of $X$.
\snl
In what follows, we will use symbols like $y$ for elements in $Y$ and then, if $a\in A$, we will write $ay$ for $y(a)$. The defining property of elements of $Y$ now reads as associativity: For all $a,a'\in A$, we have $a(a'y)=(aa')y$. Also, when $a'y$ is considered as the element $y$, acted upon by $a'$, we get
$$a(a'y)=(a'y)(a)=y(a'a)=(aa')y.$$
We see that the notations are compatible.

\inspr{1.3} Notation \rm 
Sometimes, we will use $\overline X$ for the extended module $Y$ and call it the {\it completion of $X$}.
\einspr

We refer to the Appendix where we motivate this notation and terminology in a topological context (see Proposition A.1). As we noticed already, repeating the procedure in most cases gives nothing new. This is completely in agreement with the results in the appendix.
\nl
Let us now consider some {\it basic examples}.

\inspr{1.4} Examples \rm
i) Take any set $G$ and let $A$ be the algebra $K(G)$ of complex functions with finite support on $G$ (and pointwise operations). Let $A$ act on itself by multiplication. The completed module is the space $C(G)$ of all complex functions and also the extended action is still the product. Observe that a product of any complex function with a function of finite support again has finite support.
\snl
ii) Let $A$ be any algebra and let it act on itself by multiplication on the left. The completed module consists of the right multipliers $R(A)$ of $A$, i.e.\ of 'elements' $y$ such that $ay\in A$ for all $a\in A$. The previous example is a special case of this one.
\snl
iii) Again, let $A$ be any algebra and $V$ a vector space. Let $X=A\ot V$ and let $A$ act on $X$ by multiplying the first factor from the left. So $a(a'\ot v)=(aa')\ot v$ when $a,a'\in A$ and $v\in V$. If the space $V$ is finite-dimensional, the extended module $Y$ will be $R(A)\ot V$ where $R(A)$ are the right multipliers of $A$ (as in ii). However, if $V$ is infinite-dimensional, we may get a bigger space.
\snl
iv) Suppose e.g.\ that $A=K(G)$ as in i). If $X=A\ot V$ as in iii), then $X=K(G,V)$, the space of functions from $G$ to $V$ with finite support. Then the completed module $Y$ will be $C(G,V)$, the space of all functions from $G$ to $V$. Indeed, if we multiply such a function with a function in $K(G)$, we get a function in $K(G,V)$. In general, $C(G,V)$ will be strictly bigger than $C(G)\ot V$.
\einspr

These  simple examples are important, not only for understanding the basic ideas, but also for illustrating the reason why these extended modules are important when working with algebras without identity. The last example e.g.\ shows that (algebraic) tensor products are no longer sufficient when dealing with these objects.
\snl 
We will give more, and more complicated examples later.
\nl
\it Bimodules \rm
\nl
It is clear that we can do the same thing for right $A$-modules. 
\snl
However, let us now consider the case of a bimodule. So, in what follows, $A$ will be an algebra as before, but $X$ will be a left-right $A$-bimodule. Again the left and right actions of $a\in A$ on $x\in A$ will be written as $ax$ and $xa$ respectively.
\snl
We can apply the previous basic construction, on the left module and on the right module. However, we would like to extend the bimodule structure. This is done as follows:

\inspr{1.5} Definition \rm 
Let $X$ be a (non-degenerate) $A$-bimodule. Denote by $Z$ the space of pairs $(\lambda,\rho)$ of linear maps from $A$ to $X$ satisfying
$$a\lambda(a')=\rho(a)a'$$
for all $a,a'\in A$.
\einspr

When we have such a pair of maps and if $a,a',a''\in A$, we clearly have
$$\rho(aa')a''=aa'\lambda(a'')=a\rho(a')a''$$
and because $X$ is assumed to be non-degenerate (as a right module), we find $\rho(aa')=a\rho(a')$. Similarly, because $X$ is also assumed to be non-degenerate as a left module, we get $\lambda(aa')=\lambda(a)a'$. It follows that $\rho$ belongs to the extension of the left module (as in Definition 1.1) and that $\lambda$ is an element of the extension of the right module. We will come back to this in the next item of this section.
\snl
Now, we proceed as for the basic case in the first item.
\snl
For any $x\in X$, we can define an element $(\lambda_x,\rho_x)$ in $Z$ by 
$$\lambda_x(a)=xa \qquad\qquad \text{and} \qquad\qquad \rho_x(a)=ax$$
when $a\in A$. Indeed we have because of the bimodule property,
$$a\lambda_x(a')=a(xa')=(ax)a'=\rho_x(a)a'$$
whenever $a,a'\in A$ and $x\in X$. Moreover, the map $x\mapsto (\lambda_x,\rho_x)$ from $x$ to $Z$ is injective because $X$ is assumed to be a non-degenerate module. So, we can view $X$ as sitting in $Z$.
\snl
And just as for the basic extension, also here we can extend the module action:

\inspr{1.6} Proposition \rm
For any element $z=(\lambda,\rho)$ in $Z$ and $a\in A$, we define
$$\align az&=(a\lambda(\,\cdot\,),\rho(\,\cdot\,a)) \\
         za&=(\lambda(a\,\cdot\,),\rho(\,\cdot\,)a)).
\endalign$$
This makes $Z$ into a left-right $A$-bimodule. The embedding of $X$ in $Z$, defined as above, is a bimodule map.
\snl\bf Proof\rm:
For $a,a',a''\in A$ we have
$$a'(a\lambda(a''))=aa'\lambda(a'')=\rho(a'a)a''$$
and this shows that $az\in Z$. Similarly
$$a'\lambda(aa'')=\rho(a')aa''=(\rho(a')a)a''$$
showing that also $za\in Z$.
\snl
A straightforward application of the definitions gives that 
$$(az)a'= (a\lambda(a'\,\cdot\,),\rho(\,\cdot\,a)a') =a(za')$$
for all $a,a'\in A$ and $z\in Z$. It follows that we have bimodule.
\snl
Finally, if $a\in A$ and $x\in X$ and if $z=(\lambda_x,\rho_x)$ as defined before, we get
$$\align az&=(a\lambda_x(\,\cdot\,),\rho_x(\,\cdot\,a)) \\
         za&=(\lambda_x(a\,\cdot\,),\rho_x(\,\cdot\,)a)).
\endalign$$
Because 
$$\align  a\lambda_x(a')&= axa' = \lambda_{ax}(a') \\
          \rho_x(a'a)   &= a'ax = \rho_{ax}(a')    \\
          \lambda_x(aa')&= xaa' = \lambda_{xa}(a')\\
          \rho_x(a')a   &= a'xa = \rho_{xa}(a'),
\endalign$$ 
we see that 
$$\align az & = (\lambda_{ax},\rho_{ax}) \\
         za & = (\lambda_{xa},\rho_{xa}).
\endalign$$
Therefore, the map $x\mapsto (\lambda_x,\rho_x)$ is a module map from $X$ to $Z$.
\hfill $\square$
\einspr

If $1\in A$, again we will have $Z=X$. Indeed, from the defining property, we get for any $z=(\rho,\lambda)\in Z$ 
that $\rho(1)=\lambda(1)$. And if we denote this element with $x$, we find $\rho(a)=ax$ and $\lambda(a)=xa$ for all $a$. Therefore $z=(\rho_x,\lambda_x)$.
\snl
In this case, contrary to the situation with a left or a right module, the extended module $Z$ will always we non-degenerate. Indeed, if $z=(\rho,\lambda)\in Z$ and if $az=0$ for all $a$, then we get $a\lambda(\,\cdot\,)=0$ for all $a$ and hence $\lambda=0$. Then also $\rho=0$ and so $z=0$. Similarly on the other side.
\snl
We always have $AZA\subseteq X$ and if $X$ is unital, we get equality.
\snl
As we have done before, also here we will simply write
$$az=\rho(a)\qquad\qquad\text{and}\qquad\qquad za=\lambda(a)$$
when $z=(\lambda,\rho)\in Z$ and $a\in A$. The defining relation in Definition 1.5 again reads now as associativity $(az)a'=a(za')$. One also verifies easily that this notation is compatible with the module structure defined on $Z$ in Proposition 1.6. We have e.g.\ $az=(a\lambda(\,\cdot\,),\rho(\,\cdot\,a))$ for $z=(\lambda,\rho)$ and then
$$\align a\lambda(a') &= a(za') = (az)a' \\
         \rho(a'a)    &= (a'a)z = a'(az)
\endalign$$ 
for all $a$. Similarly for $za$.
\nl
Let us now consider again some {\it examples}:

\inspr{1.7} Examples \rm
i) Let $A$ be any algebra and consider left and right multiplication. The above construction yields the multiplier algebra $M(A)$ as it was first introduced in [VD1].
\snl
ii) Again, consider any algebra $A$ and make $A\ot A$ into a bimodule by
$$\align a(a'\ot a'') &=aa'\ot a''\\
         (a'\ot a'')a &=a'\ot a''a
\endalign$$
when $a$, $a'$ and $a''$ are elements in $A$. Now, the construction will give the subalgebra of elements $x$ in $M(A\ot A)$
with the property that
$$\align (a\ot 1)x &\in A\ot A\\
         x(1\ot a) &\in A\ot A
\endalign$$
for all $a\in A$. If $(A,\co)$ is a multiplier Hopf algebra, it is a condition that $\co(A)$ belongs to this subalgebra of $M(A\ot A)$ (see [VD1]).
\snl
iii) Now, consider an algebra $A$ and a vector space $V$. Let $X=A\ot V$ and make $X$ into a bimodule in the obvious way:
$$\align a(a'\ot v) &= aa'\ot v \\
         (a'\ot v)a &= a'a\ot v
\endalign$$
where $a,a\in A$ and $v\in V$. If $A$ is a regular multiplier Hopf algebra and $\Gamma$ a coaction of $A$ on $V$, it is one of the requirements that $\Gamma$ maps $V$ into the extended module for this bimodule $X$ (see e.g.\ [De-VD-W]). We will come back to this example later. In fact, this example will appear several times in what follows, in different circumstances.
\einspr

It should be mentioned here that we can also consider the case with two algebras $A$ and $B$ and a left $A$-, right $B$-bimodule $X$. The result is again an $A$-$B$-bimodule in the obvious way.
\nl
\it Bimodules as left and right modules \rm
\nl
Let $X$ be an $A$-bimodule. Then $X$ is a left $A$-module and we can extend this to the left module, denote here by $Y_\ell$. On the other hand, $X$ is also a right module and we denote its extension with $Y_r$. If we use $Z$ to denote the extension of the bimodule as before, we have natural embeddings of $Z$ in $Y_\ell$ and of $Z$ in $Y_r$. We will now argue that in some sense, $Z$ can be considered as the intersection of $Y_\ell$ and $Y_r$.  
\snl
First observe the following. Suppose that $z\in Z$ and again, think of $z$ as a pair of maps $(\lambda,\rho)$. Then $\rho$ is completely determined by $\lambda$ and vice versa. Indeed, if $\lambda$ is given and $\rho_1, \rho_2$ are two linear maps from $A$ to $X$ satisfying
$$\align a\lambda(a')&=\rho_1(a)a' \\  
         a\lambda(a')&=\rho_2(a)a',
\endalign$$
then $\rho_1(a)a'=\rho_2(a)a'$ for all $a, a'$ in $A$ so that $\rho_1(a)=\rho_2(a)$ for all $a$ (by the assumption that $X$ is non-degenerate as a right module). So $\rho_1=\rho_2$. Similarly, when $\rho$ is given, then $\lambda$ is uniquely defined if it exists.
\snl
Therefore, we can formulate the following result:

\inspr{1.8} Proposition \rm Take an element $z$ in $Z$ and consider the associated maps $\rho, \lambda: X\to A$ given by $\rho(a)=az$ and $\lambda(a)=za$. Then we obtain an element in $Y_\ell$ and an element in $Y_r$. Conversely, if we have an element $\rho$ in $Y_\ell$ such that there exists an element $\lambda$ in $Y_r$ related by the condition $a\lambda(a')=\rho(a)a'$ for all $a,a'$, then there is a unique $z\in Z$ such that $\rho(a)=az$ and $\lambda(a)=za$ as before. Similarly, if for a given $\lambda$ there exists a $\rho$.
\einspr

This is a correct, but somewhat complicated way to say that we can view $Z$ as the intersection of $Y_\ell$ and $Y_r$.
\snl
Again, we can also consider the case of two algebras $A$ and $B$ and an $A$-$B$-bimodule.
\nl
\it More complicated and iterated cases \rm
\nl
Let us show, by considering an important example, how the situation can become quite involved when looking at more complicated cases. Moreover, because there are so many different possibilities, it seems more appropriate to consider an example to explain the techniques.
\snl
We take the case of a coaction of a multiplier Hopf algebra $A$ on a vector space $V$. Recall the following definition from [De-VD-W].

\inspr{1.9} Definition \rm
Let $(A,\co)$ be a regular multiplier Hopf algebra and assume that $V$ is a vector space. Consider the $A$-bimodule $A\ot V$ as in Example 1.7.iii and denote the extended module by $M_0(A\ot V)$. A left coaction of $A$ on $V$ is a linear map $\Gamma:V\to M_0(A\ot V)$ satisfying 
$$(\iota\ot\Gamma) \Gamma =(\co\ot\iota)\Gamma.$$
\einspr

Of course, the main difficulty is to give a precise meaning to the last equation. 
\snl
There are several possible ways to do this, but within the scope of this paper, we will give the formula a meaning by extending the maps $\iota\ot\Gamma$ and $\co\ot\iota$ to the larger space $M_0(A\ot V)$. The range of these extensions will be the space $M_0(A\ot A\ot V)$ that is obtained by applying the extension procedure to the $(A\ot A)$-bimodule $A\ot A\ot V$ (as in example 1.7.iii but with $A\ot A$ in the place of $A$). We get the following:

\iinspr{1.10} Proposition \rm
The maps $\iota\ot \Gamma$ and $\co\ot\iota$, defined on $A\ot V$, have natural extensions to maps from $M_0(A\ot V)$ to $M_0(A\ot A\ot V)$.

\snl\bf Proof\rm: 
We will in this proof write the module actions of $A$ on  $A\ot V$ and on the extended module $M_0(A\ot V)$ as 'multiplications' $(a\ot 1)y$ and $y(a\ot 1)$ when $a\in A$ and $y\in M_0(A\ot V)$. Similarly for the module actions of $A\ot A$ on $A\ot A\ot V$ and $M_0(A\ot A\ot V)$.
\snl
We begin with the map $\iota\ot \Gamma$. Take $y\in M_0(A\ot V)$ and define $z\in M_0(A\ot A\ot V)$ by
$$\align (a\ot a'\ot 1)z &= (1\ot a'\ot 1)(\iota\ot \Gamma)((a\ot 1)y) \\
         z(a\ot a'\ot 1) &= (\iota\ot\Gamma)(y(a\ot 1))(1\ot a'\ot 1).
\endalign$$
Observe that the elements $(a\ot 1)y$ and $y(a\ot 1)$ belong to $A\ot V$ so that we can apply $\iota\ot\Gamma$ to these elements. The results are elements in $A\ot M_0(A\ot V)$ and by multiplying with an element $1\ot a'\ot 1$, left or right, we get something in $A\ot A\ot V$. One easily verifies that these two linear maps from $A\ot A$ to $A\ot A\ot V$ verify the basic relation and so we get an element $z\in M_0(A\ot A\ot V)$. This gives the extension of the map $\iota\ot\Gamma$.
\snl
To define $(\co\ot\iota)(y)$ for $y\in M_0(A\ot V)$, we proceed as follows. Take $a,a'$ in $A$ and write $a\ot a'=\sum_i (p_i\ot 1)\co(q_i)$ with the $p_i$ and $q_i$ in $A$. This is possible because we have a multiplier Hopf algebra. Now define
$$(a\ot a'\ot 1)z=\sum_i (p_i\ot 1\ot 1)(\co \ot\iota)((q_i\ot 1)y)).$$
Similarly, write $a\ot a'=\sum_j\co(r_j)(s_j\ot 1)$ with $r_j$ and $s_j$ in $A$. This is possible because we have a regular multiplier Hopf algebra. Now define
$$z(a\ot a'\ot 1)=\sum_j (\co\ot\iota)(y(r_j\ot 1))(s_j\ot 1\ot 1).$$
Because $y\in M_0(A\ot V)$, we have $(q_i\ot 1)y$ and $y(r_j\ot 1)$ in $A\ot V$ and so we can apply $\co\ot\iota$. This results in elements in $M(A\ot A)\ot V$. Multiplying with $p_i\ot 1\ot 1$ from the left of $s_j\ot 1\ot 1$ from the right brings these elements down again to $A\ot A\ot V$.
\snl
So, we have defined two linear maps from $A\ot A$ to $A\ot A\ot V$. One again verifies easily that these maps satisfy the basic formula and hence, we have an element $z\in M_0(A\ot A\ot V)$. This defines the map $\co\ot\iota:M_0(A\ot V)\to M_0(A\ot A\ot V)$.

\hfill $\square$
\einspr

With these extensions, we can formulate the rule $$(\iota\ot\Gamma) \Gamma =(\co\ot\iota)\Gamma$$ needed to define coactions as in Definition 1.9.
\snl
However, more can be said about the ranges of these extensions. We will do this in the next proposition. Later, we will explain why this is important.

\iinspr{1.11} Proposition \rm When $\Gamma$ is a coaction as in Definition 1.9, we have 
$$(\iota\ot\Gamma)(M_0(A\ot V))\subseteq M_0(A\ot M_0(A\ot V).$$

\snl\bf Proof\rm: 
We have seen the spaces $M_0(A\ot V)$ and $M_0(A\ot A\ot V)$ as extensions of the bimodules $A\ot V$ and $A\ot A\ot V$ respectively. We can also look at the $A$-module $A\ot M_0(A\ot V)$ where now the last factor plays the role of $V$. We denote the extended module here with $M_0(A\ot M_0(A\ot V))$ which is in agreement with earlier conventions.
\snl
Loosely speaking, elements $z$ in $M_0(A\ot A\ot V)$ have the property that
$$\align (a\ot a'\ot 1)z &\in A\ot A\ot V \\
         z(a\ot a'\ot 1) &\in A\ot A\ot V
\endalign$$ 
for all $a,a'$, 
whereas elements $z$ in the other space $M_0(A\ot M_0(A\ot V))$ satisfy
$$\align  (a\ot 1\ot 1)z &\in A\ot M_0(A\ot V)\\
          z(a\ot 1\ot 1) &\in A\ot M_0(A\ot V)
\endalign$$
for all $a$. We see from this that we have a natural inclusion
$$M_0(A\ot M_0(A\ot V))\subseteq M_0(A\ot A\ot V).$$
If now, we take another look at the proof of the previous proposition, where we have extended the map $\iota\ot\Gamma$, we see that indeed, $(\iota\ot\Gamma)(M_0(A\ot V))\subseteq M_0(A\ot M_0(A\ot V))$.

\hfill $\square$
\einspr

In general, these two spaces are different and indeed, sometimes it may be important to know that the range of $\iota\ot\Gamma$ actually belongs to this smaller space.
\snl
Remark that there seems to be no way to obtain a similar result for the range of $\co\ot\iota$ on $M_0(A\ot V)$.
\nl
We have seen in Example 1.7.i and 1.7.ii that we can consider $A\ot A$ as an $(A\ot A)$-bimodule. This gives $M(A\ot A)$ for the extended module. But we also have other module structures (like in Example 1.7.ii) and they will yield smaller extended modules.
\snl
A similar phenomenon is encountered when considering the three-fold tensor product $A\ot A\ot A$. We get $M(A\ot A\ot A)$ when we consider the full module structure. But we also have various $(A\ot A)$-bimodule structures. Take e.g.\ multiplication of the first two factors. Let us use $M_0((A\ot A)\ot A)$ to denote the resulting extension. On the other hand, we can also consider $M_0(A\ot M_0(A\ot A))$ with constructions obtained from the obvious $A$-bimodules. We have elements $y\in M_0(A\ot A)$ when $y\in M(A\ot A)$ and
$$\align y(a\ot 1) &\in A\ot A \\ (a\ot 1)y &\in A\ot A. \endalign$$
Then elements in $M_0(A \ot M_0(A\ot A))$ are the elements $z\in M(A\ot A\ot A)$ with the property that
$$\align z(a\ot 1\ot 1) &\in A\ot M_0(A\ot A) \\
         (a\ot 1\ot 1)z &\in A\ot M_0(A\ot A). 
\endalign$$
By definition, such elements will also be in $M_0(A\ot A\ot A)$ and so we have inclusions
$$M_0(A\ot M_0(A\ot A))\subseteq M_0((A\ot A)\ot A) \subseteq M(A\ot A\ot A)$$
but in general, these inclusions are strict.
\snl
If now $(A,\Delta)$ is a regular multiplier Hopf algebra, then $(\iota\ot \Delta):A \to M_0(A\ot M_0(A\ot A))$ because
for $a,a'\in A$ we have
$$\align (a'\ot 1\ot 1)(\iota\ot\Delta)\Delta(a) &= (\iota\ot\Delta)((a'\ot 1)\Delta(a)) \\
         & \subseteq (\iota\ot\Delta)(A\ot A) = A\ot \Delta(A) \\
         & \subseteq A \ot M_0(A\ot A)
\endalign$$
and similarly on the other side. In a similar way, we will have 
$$(\Delta\ot\iota)\Delta: A\to M_0(M_0(A\ot A)\ot A)$$
where now the extensions are obtained with $A$-bimodule structures, coming from multiplication in the second factors.
\snl
This observation will turn out to be important for the use of the Sweedler notation (see the next section). 
\nl\nl

\bf 2. Coverings, mapping extensions and the Sweedler notation \rm
\nl
In the previous section, we started with an algebra $A$, with a non-degenerate product, and a left, a right, or an $A$-bimodule $X$. We showed how to get the extended module. In this section, we will moreover look at linear maps $F$ on $X$ and try to {\it extend these maps} to the extended module. For convenience, we will now assume that not only the product on $A$ is non-degenerate, but that we also have local units. Here is the precise definition.

\inspr{2.1} Definition  \rm 
We say that the algebra $A$ has {\it right local units} if, given a finite number of elements $a_1, a_2, \dots , a_n$  in $A$, there exists an element $e$ in $A$ so that $a_i e= a_i$ for all $i$. Similarly we define left local units. We say that the algebra has two-sided local units if there exist left and right local units.
\einspr

Remark that for a regular multiplier Hopf algebra, we do have left and right local units. In the case of algebraic quantum groups (regular multiplier Hopf algebras with integrals), we can even have a slightly stronger result. In this case, given a finite number of elements $(a_i)$ in $A$, there is an element $e\in A$ so that both $a_ie=a_i$ and $ea_i=a_i$ for all $i$ (see Proposition 2.6 in [Dr-VD-Z]). We will not need this stronger result.
\snl
In what follows we will assume left or/and right local units, depending on the situation. However, we should mention that the condition is not really necessary for some of the more important things we want to do. We will make various remarks about this further in this section. But it is certainly very convenient to have this condition satisfied. Also remark that the existence of local units implies that the product is non-degenerate (and that $A^2=A$).
\snl
Next, we assume that $X$ is a non-degenerate left $A$-module as before and also that $A$ has right local units. As before,  we denote the completion as defined in Definition 1.1 and Proposition 1.2 by $Y$.
\nl
{\it The covering notion: Definition and results}

\inspr{2.2} Definition \rm
Let $V$ be any vector space and let $F$ be a linear map from $X$ to $V$. We say that the variable $x$ in the expression $F(x)$ is {\it covered} from the left if there exists an element $e\in A$ such that $F(ex)=F(x)$ for all $x\in X$. 
\einspr

Of course, it is important that $e$ does not depend on $x$.
\snl
We will give plenty of examples later, but first we state and prove the {\it main} (although rather easy) {\it result}.

\inspr{2.3} Proposition \rm
Let $F: X\to V$ be a linear map from the left $A$-module $X$ to a vector space $V$ and assume that the variable is covered from the left. Then $F$ can be {\it extended} to a linear map $G:Y \to V$ satisfying $G(y)=F(ey)$ for any $y\in Y$ and any $e$ in $A$ that satisfies $F(e\,\cdot\,)=F$. This extension is {\it unique}.

\snl \bf Proof\rm:
Assume that $e_1$ and $e_2$ are elements in $A$ such that $F(e_1 x)=F(x)$ and $F(e_2 x)=F(x)$ for all $x\in X$. Choose $h\in A$ such that $e_1 h=e_1$ and $e_2 h= e_2$. Then we have for any $y\in Y$:
$$\align F(e_1y) &= F(e_1 (hy)) = F(hy) \\
          F(e_2 y)&= F(e_2(hy)) = F(hy).
\endalign$$
Therefore, we can define $G(y)=F(ey)$ where $e$ is any element in $A$ such that $F(e\,\cdot\,)=F$.
\snl
It is clear that $G$ is linear and also the uniqueness follows from the argument given above.
\hfill $\square$
\einspr

Also the case of a linear map on a {\it right $A$-module} can be considered and of course, we have similar definitions and results. In this case, we assume that $A$ has {\it left local units}.
\snl
The extension can also be characterized in a topological way (see Proposition A.2 in the appendix). For obvious reasons, we will often denote the extended map with the same symbol.
\nl
Let us now look at the case of a {\it bimodule}. We will work with two algebras $A$ and $B$ and a non-degenerate $A$-$B$-bimodule $X$ with $A$ acting on the left and $B$ on the right. We assume that $A$ has right local units and that $B$ has left local units. As before, we denote the completion, as defined in Definition 1.5 and Proposition 1.6, by $Z$.  
\snl
Again, suppose that $F:X \to V$ is a linear map from the space $X$ to any vector space $V$. The following gives the appropriate covering notion.

\inspr{2.4} Definition \rm
We say that the variable $x$ in the expression $F(x)$ is covered if it is covered both from the left and from the right.
\einspr              

So, we mean that there exists $e\in A$ and $f\in B$ so that $F(ex)=F(x)$ and $F(xf)=F(x)$ for all $x$. Observe that this is a {\it stronger} assumption than requiring the existence of elements $e,f$ such that $F(exf)=F(x)$ for all $x$. This would be a too weak assumption for our purposes.
\snl
Because we have a covering left and right, we can apply Proposition 2.3 and extend $F$ to $Y_\ell$ and $Y_r$, the completions of the left and the right modules respectively. The following implies that these extensions coincide on the 'intersection' $Z$ of $Y_\ell$ and $Y_r$ (remember the discussion in Section 1, see Proposition 1.8).

\inspr{2.5} Proposition \rm
Let $F:X\to V$ be a linear map from the $A$-$B$-module $X$ to the vector space $V$. Assume that the variable is covered. Then there exists a {\it unique extension} of $F$ to a linear map $G:Y\to V$ such that 
$$G(y)=F(ey)=F(yf)$$
for all $y$ in $Y$ and for any pair of elements $e\in A$ and $f\in B$ that satisfy $F(e\,\cdot\,)=F$ and $F(\,\cdot\,f)=F$.

\snl \bf Proof\rm: 
If $e$ and $f$ are as in the formulation above, then for all $y$ in $Y$ we have
$$F(ey)=F(eyf)=F(yf).$$
This, in combination with the arguments in the proof of Proposition 2.3, will yield the result.
\hfill $\square$
\einspr

Again, this extension is unique and we can characterize it by the requirement that $F(y)=F(ey)$ when $e$ is such that $F(e\,\cdot\,)=F$ or that $F(y=F(yf)$ when $f$ satisfies $F(\,\cdot\,f)=F$. Moreover, also in this case, we will often use the same symbol for the extended map. There is also a {\it topological interpretation} and then, the extension is the unique {\it continuous} extension (see Proposition A.3 in the appendix). 
\nl
These results are of not very deep. The examples are more interesting. And before we continue, let us look at two important examples to illustrate this technique. One is about extending an {\it action} of $A$ to its multiplier algebra $M(A)$. The other one is about extending a {\it pairing} of multiplier Hopf algebras. We will give more examples later in this section and in the next one.

\inspr{2.6} Examples \rm
i) Take a left $A$-module $X$ and assume that it is unital. This means that $AX=X$. If $A$  has left local units, this implies that for any $x_1, x_2, \dots , x_n$ there is an element $e\in A$ such that $ex_j=x_j$ for all $j$. Now consider $A$ as an $A$-bimodule. Its completion is $M(A)$. Given $x\in X$, define $F:A\to X$ by $F(a)=ax$. The variable is covered from the right because there exits an element $e\in A$ such that $ex=x$. When we apply Proposition 2.3 (for the right module), we find an extension $G$. And if we restrict this extension again to $M(A)$, we see that we can define $mx$ for $m\in M(A)$ and $x\in X$ such that $mx=mex$ whenever $e\in A$ satisfies $ex=x$. This is how a unital left $A$-module is extended to a unital left $M(A)$-module (see Proposition 3.3 in [Dr-VD-Z]).
\snl
ii) Let $\langle\,\cdot\, , \,\cdot\,\rangle : A \times B \to \Bbb C$ be pairing of regular multiplier Hopf algebras (in the sense of [Dr-VD]). Consider $A$ as an $A$-bimodule. Fix $b\in B$ and define $F:A\to \Bbb C$ by $F(a)=\langle a,b \rangle$. We know that the left and the right actions of $A$ on $B$, induced by the pairing, are unital. If we denote these actions by $a\tr b$ and $b\tl a$ respectively, we find,  elements $e,f\in A$ such that $f\tr b=b$ and $b\tl e=b$. Then
$$\align F(ea)&=\langle ea,b\rangle = \langle a, b\tl e\rangle = \langle a,b \rangle \\
         F(af)&=\langle af,b\rangle = \langle a, f\tr b\rangle = \langle a,b \rangle 
\endalign$$
for all $a$ and we see that $F(e\,\cdot\,)=F(\,\cdot\,f)=F$. So, the variable is covered (through the pairing by $b$). Then applying Proposition 2.5, we can extend the map $F$ to a linear map on $M(A)$. We use $\langle m,b\rangle$ to denote the value of the extension in the element $m\in M(A)$. We find that 
$$\align \langle m, a\tr b\rangle &= \langle ma,b \rangle \\
         \langle m, b\tl a\rangle &= \langle am,b \rangle
\endalign$$
for all $a\in A$, $b\in B$ and $m\in M(A)$. This is how  a pairing between $A$ and $B$ is extended to a bilinear map on $M(A) \times B$ (see e.g.\ [De]). Similarly, it can be extended to $A\times M(B)$. Remark however that in general, it is not possible to extend to $M(A) \times M(B)$.
\einspr  

We are now ready for the following, important item.
\nl
\it The Sweedler notation \rm
\nl
Essentially, we have to look at two different cases. We have the usual {\it coproducts} and there are the {\it coactions}. Let us first consider the coproducts. For this, we assume that $(A,\Delta)$ is a multiplier Hopf algebra.

\inspr{2.7} Notation \rm
Consider $A\ot A$ as a left $A$-module where the action is given by left multiplication in the first factor. Let $V$ be a vector space and let $F:A\times A\to V$ be a bilinear map. Consider $F$ also as a linear map on $A\ot A$ and assume that the variable is covered from the left (as in Definition 2.2). Then we can extend $F$ to the completed module. Because this contains $\Delta(A)$, we can define $F(\Delta(a))$ for $a\in A$. We write
$\sum_{(a)}F(a_{(1)},a_{(2)})$
to denote this element in $V$. \hfill$ \square$
\einspr

So, by definition we have 
$$ \sum_{(a)}F(a_{(1)},a_{(2)})=\sum_i F(p_i,q_i)$$
where $\sum_i p_i\ot q_i=(e\ot 1)\Delta(a)$ and where $e\in A$ is chosen so that $F(e\,\cdot\,)=F$. We say that, in the expression $\sum_{(a)}F(a_{(1)},a_{(2)})$, the factor $a_{(1)}$ is covered from the left.
\snl
There are clearly many other possibilities, but let us first illustrate the above with two examples.

\inspr{2.8} Example \rm
i) As before, let $(A,\Delta)$ be a multiplier Hopf algebra. Fix $b\in A$ and define $F:A\times A\to A\ot A$ by $F(p,q)=bp\ot q$. Clearly $p$ is covered from the left (take $e$ such that $be=b$). So we can write $\sum_{(a)}ba_{(1)}\ot a_{(2)}$ (which is of course $(b\ot 1)\Delta(a)$ here). 
\snl
ii) Again, let $A$ be a regular multiplier Hopf algebra $A$ and let $R$ be a right $A$-module algebra (as in Definition 4.1 of [Dr-VD-Z]). So $R$ is a unital right $A$-module and if we denote the action of $a\in A$ on $x\in R$ by $x\tl a$, we also require that 
$$(xx')\tl a=\sum_{(a)} (x\tl a_{(1)})(x'\tl a_{(2)})$$
whenever $x,x'\in R$ and $a\in A$, using the Sweedler notation. 
\snl Let us explain this in more detail. We are taking fixed elements $x,x'\in R$ and we consider the map $F:A\ot A\to R$, defined by
$$F(a\ot a')=(x\tl a)(x'\tl a').$$
In order to get the appropriate extension of this map $F$, we consider $A\ot A$ as left $A$-module and write the module action of $a$ as multiplication with $a\ot 1$. Choose an element $e\in A$ such that $x\tl e=x$. Then $F((e\ot 1)c)=F(c)$ for all $c\in A\ot A$. It follows that the variable $c$ is covered from the left and we can apply Proposition 2.3 to extend $F$. If we restrict again to $M(A\ot A)$, we see that we can define $G:M(A\ot A)\to R$ with the property that 
$$G(m)=\sum_i (x\tl p_i)(x'\tl q_i)$$
where $e$  is as before and where $\sum_i p_i\ot q_i= (e\ot 1)m$.  We can now write 
$$G(\co(a))=\sum_{(a)} (x\tl a_{(1)})(x'\tl a_{(2)})$$
as is explained in 2.7. We see from this example why it is handy to be able to use the Sweedler notation here.
\hfill $\square$
\einspr

In the previous example, we could have considered $A\ot A$ as a left $A$ module with multiplication from the left in the second factor. Then we would have written
$$\sum_{(a)} (x\tl a_{(1)})(x'\tl a_{(2)})=\sum_i (x\tl p_i)(x'\tl q_i)$$
where now 
$\sum_i p_i\ot q_i= (1\ot e)\Delta(a)$. In this case, the factor $a_{(2)}$ is covered from the left. 
\snl
Another possibility is to consider  $A\ot A$ as left $(A\ot A)$-module. Now we choose elements $e,e'\in A$ such that $x\tl e=x$ and $x'\tl e'=x'$. Then $F((e\ot e')c)=F(c)$ for all $c\in A\ot A$. It follows that the variable $c$ is covered from the left and we can again apply Proposition 2.3 to extend $F$. If we restrict again to $M(A\ot A)$, we see that we can define $G:M(A\ot A)\to R$ with the property that 
$$G(m)=\sum_i (x\tl p_i)(x'\tl q_i)$$
where $e$ and $e'$ are as before and where $\sum_i p_i\ot q_i= (e\ot e')m$.  We can now write 
$$G(\co(a))=\sum_{(a)} (x\tl a_{(1)})(x'\tl a_{(2)}).$$
\snl
From this point of view, we see that the result in Proposition 2.5 will give that it does not matter for the final result if we take $A\ot A$ as a left $A$-module, either by considering left multiplication in the first factor, or in the second factor. 
\snl
This case is typical. In general, one, or a stronger covering, is enough to {\it define} an expression, but most of the time, different coverings are needed for {\it proving} equations. Above, the stronger covering e.g.\ with the two elements $e,e'$ is sufficient to give a meaning to the expression. It also allows to define the expression for more cases. However, the other coverings, with only one element ($e$ or $e'$) will be needed in order to {\it prove} equalities involving these expressions. We will illustrate this in example 2.10.
\nl
Still related with the notion of a module algebra is the following example about the {\it twist} map.

\inspr{2.9} Example \rm
i) As in the previous example, let $A$ be a regular multiplier Hopf algebra and $R$ a right $A$-module algebra. Consider the twist map $T:R\ot A \to A\ot R$ given by 
$$T(x\ot a)=\sum_{(a)}a_{(1)} \ot (x\tl a_{(2)})$$
where we use the Sweedler notation. To give a meaning to the right hand side, we fix $x\in R$ and consider the linear map $F:A\ot A \to A\ot R$ given by
$$F(p\ot q)=p\ot (x\tl q).$$
We look at the left $A$-module $A\ot A$ where the action is given by multiplication of the second factor from the left. We see that the variable is covered through the action by the element $x$. Then we can extend $F$ to the extended module. As $\co(a)$ belongs to this extended module for all $a\in A$, we can define the map $T:R\ot A\to A\ot R$ by
$R(x\ot a)=\sum_i p_i \ot (x\tl q_i)$
where $\sum_i p_i\ot q_i=(1\ot e)\co(a)$ and where $e$ is any element in $A$ satisfying $x\tl e=x$.
\snl
ii) Let us also look at the smash product $A\# R$. Using the Sweedler notation, the product is written as 
$$(a\# x)(a'\# x')=\sum_{(a')}a a'_{(1)} \# (x\tl a'_{(2)})x'$$
when $a,a'\in A$ and $x,x'\in R$. In order to apply our method here, we fix elements $a\in A$ and $x,x'\in R$ and consider the map  $F:A\ot A\to A\ot R$ given by
$$F(p\ot q)=ap \ot (x\tl q)x'.$$
Now, we view $A\ot A$ as an $A$-bimodule where $A$ acts by multiplication from the left, both in the first and the second factor (as in Example 2.8.ii). Then, the variable is covered in both places. For the multiplication in the first variable, we use an element $e$ such that $ae=a$. For the multiplication in the second variable, we use an element $f$ such that $x\tl f=x$. Now, we apply Proposition 2.5. We get two possible expressions. On the one hand, we can write
$$(a\# x)(a'\# x')=\sum_i a p_i \# (x\tl q_i)x'$$
where $\sum_i p_i\ot q_i=(e\ot 1)\co(a)$ and $ae=a$. On the other hand, we have
$$(a\# x)(a'\# x')=\sum_j a r_j \# (x\tl s_j)x'$$
where now $\sum_j r_j\ot s_j=(1\ot f)\co(a)$ and $x\tl f=x$.
\einspr

The idea should be clear and the reader should now be able to apply the technique in many other situations. However, let us consider one more example where now, we look at {\it coactions}. 
\snl
We consider the case of a right coaction here (as the notations are somewhat easier). So, let $(A,\co)$ be a regular multiplier Hopf algebra. Assume that $V$ is a vector space and that $\Gamma$ is a right coaction of $A$ on $V$ (cf.\ Definition 1.9). So, $\Gamma$ is an injective linear map from $V$ to $M_0(V\ot A)$, the completion of the left $A$-module $V\ot A$ where the action is multiplication, left and right, of the second factor. We also have the condition 
$$(\Gamma\ot\iota)\Gamma = (\iota\ot\co)\Gamma$$
(as explained in Propositions 1.10 and 1.11 for a left coaction). 

\iinspr{2.10} Example \rm
i) Apply the counit $\varepsilon$ of $A$ on the last factor in the coassociativity formula. Using the Sweedler notation, we have for $a\in A$:
$$\align \Gamma(\sum_{(v)}\varepsilon(v_{(1)})\Gamma(v_{(0)})) 
     &= \sum_{(v)}v_{(0)}\ot (\iota\ot\varepsilon)\Delta(v_{(1)}) \\
     &= \sum_{(v)}v_{(0)}\ot v_{(1)}=\Gamma(v).
\endalign$$
From the injectivity of $\Gamma$, we find that $(\iota\ot\varepsilon)\Gamma(v)=v$ for all $v\in V$. In the first expression, the variable $v_{(1)}$ is covered as we apply $\varepsilon$ (by multiplying with $e\in A$ such that $\varepsilon(e)=1$). If we multiply with an extra element of $A$ in the second factor, we also get $v_{(1)}$ covered in the second expression.
\snl
ii) Now, consider the cotwist map $T$ from $A\ot V$ to $V\ot A$, defined by 
$$T(a\ot v)=(1\ot a)\Gamma(v).$$
Because $\Gamma(v)\in M_0(V\ot A)$ for all $v\in V$ by definition, this map is well-defined. 
\snl
The cotwist map is known to be bijective and we will now consider the obvious candidate for the inverse. Using the Sweedler notation, it is given by 
$$v\ot a \mapsto \sum_{(v)} aS^{-1}(v_{(1)})\ot v_{(0)}$$
where $\Gamma(v)= \sum_{(v)} v_{(0)}\ot v_{(1)}$. To define this map properly, using our theory, we fix the element $a\in A$ and define the linear map $F:A\ot V \to A\ot V$ by
$$F(w\ot p)=aS^{-1}(p)\ot w.$$
If we let $A$ act on $V\ot A$ by multiplying the second factor from the right, we see that the variable is covered by $a$. Indeed, if $e\in A$ is chosen so that $ae=a$, then $aS^{-1}(p)=aS^{-1}(pS(e))$ and $F(\,\cdot\,S(e))=F$. Therefore, we can extend $F$ to the completed module. We get the map
$$v\ot a \mapsto \sum_i aS^{-1}(p_i)\ot w_i$$
where $\sum_i w_i\ot p_i = \Gamma(v)(1\ot S(e))$ and $ae=a$.
\snl
iii) Next, let us try to argue that the map defined in i) is indeed the inverse of the cotwist map $T$. Using the Sweedler notation, the formal argument is like this.
\snl
Start with $a\in A$  and $v\in V$ and consider 
$$T(a\ot v)=\sum_{(v)} v_{(0)} \ot av_{(1)}.$$
Then apply the candidate for the inverse. The result is
$$\sum_{(v)} av_{(2)}S^{-1}(v_{(1)}) \ot v_{(0)}=\sum_{(v)} \varepsilon(v_{(1)})a \ot v_{(0)}=a\ot v.$$
We use that $(\iota\ot\Delta)\Gamma=(\Gamma\ot\iota)\Gamma$, that $\sum_{(p)}p_{(2)}S^{-1}(p_{(1)})=\varepsilon(p)1$ for $p\in A$ and that $(\iota\ot\varepsilon)\Gamma(v)=v$ for $v\in V$.
\snl
Suppose now that we want to make this argument correct by applying the right coverings. In the first place, we observe that $v_{(1)}$ is covered by $a$ in the expression for $T(a\ot v)$ and so we can apply the candidate for the inverse map. In the next expression, we have a successive covering: First $v_{(2)}$ is covered from the left by $a$  and then, $v_{(1)}$ is covered (from the right because of the presence of the map $S^{-1}$). The next step is a bit more tricky. We need to use the equality $\sum_{(p)}p_{(2)}S^{-1}(p_{(1)})=\varepsilon(p)1$ for $p=v_{(1)}$. A safe way to do this would be to use an extra covering, say by multiplying with $a'$ from the right in the first factor. 
\einspr

We could have given a more direct definition of the inverse map by
$$v\ot a\mapsto \sigma (\iota\ot S^{-1})(\Gamma(v)(1\ot S(a)))$$
where $\sigma$ is used for the flip map on $V\ot A$. In the example, we have used a different way in order to illustrate the theory. This point of view will also be useful later, in Section 3, where we will come back to this case.
\nl
Very similar is the equation
$$\sum_{(a)}a'S(a_{(1)})a_{(2)} \ot a_{(3)}=a'\ot a$$
for $a,a'\in A$. It is used to obtain the inverse of the map $a'\ot a\mapsto (a'\ot 1)\Delta(a)$ on $A\ot A$. This is again a case where we have a successive (or iterated) covering. We will focus on this aspect in the next section where we begin with an example related to the above case.
\nl\nl

\bf 3. More complicated cases and more examples \rm
\nl
In the first two sections of this paper, we have formulated the basic ingredients together with the basic results and we have used relatively simple examples to illustrate these. In this section however, we will look at more complicated, in particular iterated cases and discuss also more complicated examples. We will see e.g.\ how the idea of {\it successive} (or {\it iterated}) coverings will be useful in some situations. And because there are so many different possible cases that are all alike, we will do this only by means of examples. 
\snl
Most of the more complicate examples are taken from the recent paper on the {\it bicrossproduct of multiplier Hopf algebras} (see [De-VD-W]). However, let us start with the following simplest example to illustrate the case of an iterated covering.

\inspr{3.1} Example \rm 
i) Let $A$ be a regular multiplier Hopf algebra and consider the formula
$$\sum_{(a)}a_{(1)}S(a_{(2)})a_{(3)}=a$$
for $a\in A$. When we want to cover the expression in the left hand side of this formula, we can do this if we multiply left and right with elements of $A$. This means that we look at the left $(A\ot A^{\text{op}})$-module $A\ot A\ot A$ with $A$ acting by multiplication on the left in the first factor and with $A^{\text{op}}$ acting by multiplication on the right in the third factor. Then, we take $a',a''\in A$ and look at the linear map $F:A\ot A\ot A\to A$ given by 
$$F(p\ot q\ot r)=a'pS(q)ra''.$$ It is clear that the variable is covered and so we can extend the map to the completed module. Because obviously, $(\iota\ot\Delta)\Delta(a)$ belongs to this completed module, we can write $\sum_{(a)}a'a_{(1)}S(a_{(2)})a_{(3)}a''$ as an element of $A$ for all $a,a',a''\in A$.
\snl
ii) It is more or less clear that we do not really need that much covering. It should be enough when we only multiply on one side, say on the left. Indeed, as $(a'\ot 1)\Delta(a)\in A\ot A$, we get that $\sum_{(a)}a'a_{(1)}S(a_{(2)})a_{(3)}$ is a linear combination of elements of the form $\sum_{(q)}a'pS(q_{(1)})q_{(2)}$ with $p,q\in A$ and then we see that $q_{(1)}$ will be covered from the right by the element $p$.
\snl
iii) In order to make the reasoning above precise and to put it in the general framework, we first have to consider $A\ot A$ as a left $A^{\text{op}}$-module by letting $A$ multiply on the right in the first factor. Denote the completed module by $M_0(A\ot A)$. Intuitively speaking, this consists of 'elements' $y$ with the property that $y(a\ot 1)\in A\ot A$ for all $a\in A$.  Then consider $A\ot M_0(A\ot A)$ as a left $A$-module with the action given by multiplication on the left in the first factor. Denote the completed module by $M_1(A\ot A\ot A)$. Again intuitively speaking, the 'elements' $z$ in this completion will satisfy $(a\ot 1\ot 1)z\in A\ot M_0(A\ot A)$.  
Now, we look at the map 
$$p\ot q\ot r \mapsto a'pS(q)r$$ on 
$A\ot A\ot A$ with $a'$ given in $A$. When we fix $p$, we get a map on $A\ot A$. The variable is covered and we can extend it to the completion $M_0(A\ot A)$. This results in a map from $A\ot M_0(A\ot A)$ to $A$. Again the variable will be covered and we can extend it to $M_1(A\ot A\ot A)$. Finally observe that $(\iota\ot \Delta)\Delta(a)\in M_1(A\ot A\ot A)$ for all $a\in A$ because, for any other $a'$ in $A$, we get
$$(a'\ot 1\ot 1)(\iota\ot \Delta)\Delta(a)=(\iota\ot \Delta)((a'\ot 1)\Delta(a))\subseteq A\ot \Delta(A)$$
and $\Delta(A)\in M_0(A\ot A)$. For all these reasons, it makes sense to define the element 
$\sum_{(a)}a'a_{(1)}S(a_{(2)})a_{(3)}$ in $A$ for all $a,a'\in A$. \hfill $\square$
\einspr

We could also consider $A\ot A\ot A$ as a left $(A\ot A^{\text{op}})$-module    
by letting $A$ and $A^{\text{op}}$ act as above. The completed module $M_0(A\ot A\ot A)$ will consist of 'elements' $z$ satisfying $(a\ot 1\ot 1)z(1\ot a'\ot 1)\in A\ot A\ot A$ for all $a,a'\in A$. This will (in general) be strictly bigger than $M_1(A\ot A\ot A)$. If we do this however, the variable will not necessarily be covered for the map $p\ot q\ot r \mapsto a'pS(q)r$ under consideration and therefore, the method can not be used to extend the map to this possibly bigger completion.
\snl
Observe that we already had a case of successive covering in the last example of the previous section, Example 2.10.iii. Also remark that we have not given a rigorous proof of the equality in Example 3.1. We just have explained how to give a meaning to the expression $\sum_{(a)}a_{(1)}S(a_{(2)})a_{(3)}$. In order to show that actually, this is equal to $a$, we can use more coverings than needed to define the expression. In this case, we can use ideas as above and as in the last example of the previous section. In fact, the equality suggests that the expression needs no extra covering at all. This can be seen e.g.\ if we multiply with $\varepsilon(a')$ (with $a'$ taken so that this is $1$), inserting it between $a_{(1)}$ and $S(a_{(2)})$.
\nl
Let us now look at an example with reference to the definition of the {\it smash coproduct} $\Delta_\#$ on the smash product $AB$ of multiplier Hopf algebras $A$ and $B$, considered in Proposition 3.2. in [De-VD-W].   

\inspr{3.2} Example \rm
i) Recall that $A$ and $B$ are regular multiplier Hopf algebras with a left coaction $\Gamma:A\to M(B\ot A)$, making $A$ into a left $B$-comodule coalgebra. We have defined $\Delta_\#(a)=\sum_{(a)}(a_{(1)}\ot 1)\Gamma(a_{(2)})$ in $M(AB\ot AB)$ for $a\in A$. The algebra $AB$ is obtained as the smash product constructed with a right action of $A$ on $B$ making $B$ into an $A$-module algebra.
\snl
If we multiply with an element $a'\in A$ from the left in the first factor, then $a_{(1)}$ will be covered. We get an element in $(A\ot 1)\Gamma(A)$. Next, if we multiply with an element $b\in B$, again from the left in the first factor, and if we use that $BA=AB$, we end up with an element in $(AB\ot 1)\Gamma(A)$. So, the first leg of $\Gamma$ will also be covered. This  argument is given in the first part of the proof of Proposition 3.2 of [De-VD-W].
\snl
ii) Let us now see, more precisely, how this is an iterated covering. We take $a'\in A$ and $b\in B$ and we consider the linear map $F$ from $A\ot B\ot A$ to $AB\ot A$ given by 
$$F(p\ot q\ot r)=ba'pq\ot r.$$
First we fix $p\in A$. We consider $B\ot A$ as a left $B$-module (with multiplication from the left in the first factor) and we consider the completed module $M_0(B\ot A)$. For the restricted map from $B\ot A$ to $AB\ot A$, given by
$$q\ot r \mapsto F(p\ot q\ot r)=ba'pq\ot r,$$
the variable is covered (as $ba'p\in AB$) and we can extend the map to the completed module $M_0(B\ot A)$. As $\Gamma(a)$ belongs to this completion, we can define
$$F(p\ot \Gamma(a))=(ba'p\ot 1)\Gamma(a)$$
in $AB\ot A$ for all $a\in A$. This is of course nothing new. 
\snl
For the next step, we consider the left $A$-module $A\ot M_0(B\ot A)$ with again $A$ acting from the left in the first variable. Denote the completed module with $M_1(A\ot B\ot A)$. For the map from $A\ot M_0(B\ot A)$ to $AB\ot A$, given by
$$p\ot y \mapsto F(p\ot y)=(ba'p\ot 1)y,$$ 
the variable is covered (by $a'$). So, we can extend it to the completed module. Because $(\iota\ot \Gamma)\Delta(a)$ belongs to this module, we can define the element
$$\sum_{(a)}  (ba'a_{(1)}\ot 1)\Gamma(a_{(2)})$$
in $M(AB\ot AB)$ for $a\in A$. \hfill $\square$
\einspr

Of course, for the above example, we could simply look immediately at the left $A$-module $A\ot A$ (with multiplication from the left in the first factor) and consider the function from $A\ot A$ to $AB\ot A$, given by $p\ot q\mapsto (ba'p\ot 1)\Gamma(q)$. Then, the covering of the first leg of $\Gamma$ is already taken into account. 
\snl
Another example of such an iterated covering is found in the definition of the map $P$ in Proposition 3.3 of [De-VD-W]. We will look a bit closer at the definition of this map and some of its properties in Example 3.5 later in this section.
\nl
Another case of successive coverings is encountered often when using that a coproduct is a homomorphism. Here is a more complicated  example of this situation.

\inspr{3.3} Example \rm
i) We consider the proof of Lemma 3.7 of [De-VD-W]. We are in the situation as in the previous example and we also use the coproduct $\Delta_\#$ on the smash product $AB$. Also recall the definition of the right action of $AB$ on $B$,  defined by $y\bullet(ab)=(y\tl a)b$ where $a\in A$, \, $b,y\in B$ and where $\tl$ is used for the original right action of $A$ on $B$. This action is then combined with right multiplication of $A$ on itself to yield a right action of $AB\ot A$ on $B\ot A$ denoted with the same symbol and given by 
$$(y\ot x)\bullet (c\ot a)=(y\bullet c)\ot xa$$
when $a,x\in A$, $y\in B$ and $c\in AB$. See Notation 3.4 in [De-VD-W].
\snl
In the proof of that lemma, we have the expression
$$\sum_{(a)(a')} (xa_{(1)}a'_{(1)}\ot x')((y\ot 1)\bullet\Delta_\#(a_{(2)})\bullet\Delta_\#(a'_{(2)}))$$
where $a,a',x,x'\in A$ and $y\in B$. First, $a_{(1)}$ is covered by $x$ and then $a'_{(1)}$ will also be covered. The result is a linear combination of elements of the form
$$(p\ot 1)((y\ot x')\bullet\Delta_\#(q)\bullet\Delta_\#(q'))$$
with $p,q,q'\in A$.
\snl
ii) What about this iterated covering? We can look at it in the following way. Consider $A\ot A$ as a left $A$-module (with multiplication on the left in the first factor as usual) and the map $F:A\ot A \to AB\ot A$ given by
$$F(p'\ot q')= (xpp'\ot 1)((y\ot x')\bullet\Delta_\#(q)\bullet\Delta_\#(q'))$$
with $p,q\in A$. The variable is covered by $p$ and so we can extend the map to $\Delta_A(a')$. Next, consider the map $G:A\ot A\to AB\ot A$ given by 
$$G(p\ot q)= \sum_{(a')}(xpa'_{(1)}\ot 1)((y\ot x')\bullet\Delta_\#(q)\bullet\Delta_\#(a'_{(2)}))$$
with the same module structure. Again, the variable will be covered by $x$ and so we can extend the map to $\Delta_A(a)$. This yields the expression we are looking at. \hfill $\square$
\einspr

In the next example, we treat a rather complicated case of iterated coverings (related with Proposition 3.8 of [De-VD-W]).

\inspr{3.4} Example \rm 
i) Consider the formula
$$\Gamma(aa')=\sum_{(a)(a')}((a_{(-1)}\tl a'_{(1)})\ot a_{(0)})\Gamma(a'_{(2)})$$
for  $a,a'\in A$. We use the Sweedler notation for the coproduct $\Delta_A(a')$ and the coaction $\Gamma(a)$ with $a,a'\in A$.
If we multiply with an element $b$ of $B$ on the left in the first factor, the left hand side will be covered and will yield an element in $B\ot A$. So, we expect that the same will be true for the right hand side. However, this is certainly not immediately clear.
\snl
To see this is the case, observe that for all $p\in B$ and $q\in A$ we have
$$b(p\tl q)=\sum_{(q)}((b\tl S(q_{(1)}))p)\tl q_{(2)}$$
(using the $A$-module algebra property). The element $q_{(1)}$ is covered from the right by the element $b$ through the action. Then, when we look at the expression
$$\sum_{(a)(a')}(((b\tl S(a'_{(1)}))a_{(-1)})\tl a'_{(2)}\ot a_{(0)})\Gamma(a'_{(3)}),$$
we get first that $a'_{(1)}$ is covered by $b$ through the action, next $a_{(-1)}$ will be covered through multiplication with an element of $B$, then $a'_{(2)}$ will be covered, again through the action, and finally, the first leg of $\Gamma(a'_{(3)})$ is covered by multiplication by an element of $B$. So we see that in the original formula, everything seems to be  well covered when we multiply on the left in the first factor with an element of $B$. This results is an element in $B\ot A$.
\snl
ii) Let us now consider this argument in greater detail. We start with the linear map $F:B\ot A\to B\ot A$, given by
$$F(q\ot r)=\sum_{(a')}(b(q\tl a'_{(1)})\ot r)\Gamma(a'_{(2)})$$
where $a'\in A$ and $b\in B$ are given. Remark that $a'_{(1)}$ is covered through the action by $q$ and that also the first leg of $\Gamma(a'_{(2)})$ is covered. So, we really get a well-defined map from $B\ot A$ to $B\ot A$. We have
$$\sum_{(a')}b(q\tl a'_{(1)})\ot a'_{(2)}=
\sum_{(a')}((b\tl S(a'_{(1)}))q)\tl a'_{(2)}\ot a'_{(3)}.$$
In the expression on the left, we have that $a'_{(1)}$ is covered by $q$. On the right, $a'_{(1)}$ is covered by $b$ and then also $a'_{(2)}$ will be covered. In fact, take $e\in A$ so that $b\tl e=b$. Write
$$\sum_{(a')}eS(a'_{(1)})\ot a'_{(2)}=\sum_i p_i\ot q_i$$
in $A\ot A$. Then 
$$\sum_{(a')}b(q\tl a'_{(1)})\ot a'_{(2)}=
\sum_{i,(q_i')}((b\tl p_i)q)\tl q_{i(1)}\ot q_{i(2)}.$$
Next, choose $f\in B$ so that $(b\tl p_i)f=b\tl p_i$ for all $i$. We find that
$$\sum_{(a')}b(q\tl a'_{(1)})\ot a'_{(2)}=
\sum_{(a')}b((fq)\tl a'_{(1)})\ot a'_{(2)}.$$
Remark that $f$ depends on $a$ and on $b$ but not on $q$. It follows that the variable of $F$, defined above, is covered and that we can extend $F$ and define $F(\Gamma(a))$ for $a\in A$. Using the conventions about the Sweedler notation, we can write it as
$$\sum_{(a)(a')}(b(a_{(-1)}\tl a'_{(1)})\ot a_{(0)})\Gamma(a'_{(2)}).$$
\hfill $\square$
\einspr

The above example certainly illustrates how subtle these iterated coverings can be. In the last expression, $a_{(-1)}$  appears to be covered by $b$ and $a'$ but this is certainly not easy to see when simply examining the formula.
\snl
Now, as promised earlier, we look at an example related with formulas involving the map $P$ as defined in Proposition 3.3 of [De-VD-W]. The situation is very similar to the one encountered in the previous example. The result however is quite handy so that using the map $P$ becomes in a sense more easy than the twist and cotwist maps $R$ and $T$, or rather their 'opposite' counterparts $R^{\text{op}}$ and $T^{\text{op}}$ that compose to $P$ (see a remark after Proposition 3.3 in [De-VD-W]).

\inspr{3.5} Example \rm
i) So, with the ingredients as in the previous examples, consider the map $P:B\ot A\to  B\ot A$ given by
$$P(b\ot a)=\sum_{(a)}((b\tl a_{(1)})\ot 1)\Gamma(a_{(2)}).$$
 Take $a\in A$ and $b'\in B$ and define $F:B\to B\ot A$ by
$$F(q)=(b'\ot 1)P(q\ot a)=\sum_{(a)}(b'(q\tl a_{(1)})\ot 1)\Gamma(a_{(2)}).$$
We have seen in the previous example that there exists an element $f\in B$ such that $F(q)=F(fq)$ for all $q\in B$.
\snl
ii) Next consider the map $G:B\ot B\to B\ot B\ot A$ given by
$$G(q_1\ot q_2)=(b'\ot 1\ot 1)P_{23}P_{13}(q_1\ot q_2\ot a)$$
where $a\in A$ and $b'\in B$. Because of the result in i), we find an element $f\in B$ so that
$$G(q_1\ot q_2)=G(fq_1\ot q_2)$$
for all $q_1,q_2\in B$. Therefore, we can extend $G$ and define it on $\Delta_B(b)$ for any $b\in B$. This yields the expression
$$(b'\ot 1\ot 1)P_{23}P_{13}(\Delta_B(b)\ot a)$$
needed in Assumption 3.9 of [De-VD-W]. Similarly, we can define
$$(1\ot b'\ot 1)P_{23}P_{13}(\Delta_B(b)\ot a).$$
Observe that the other expression in this assumption, namely $(\Delta_B\ot\iota_A)P(b\ot a)$, is also covered when multiplying from the left with an element of $B$, either in the first or in the second factor.
\snl
iii) We can also take a quick look at the proof of Proposition 3.10 of [De-VD-W]. If we multiply all equations with an element $b\in B$ from the left in the first factor, we arrive at elements in $B\ot A$ or $B\ot B\ot A$ and we can safely apply $\iota_B\ot \varepsilon_A$ or $\iota_B\ot\iota_B\ot\varepsilon_A$ respectively.
\hfill $\square$
\einspr

So far several examples with more complicated situations, mostly taken from the paper [De-VD-W]. The idea should be clear from these examples. And the reader should now be able to understand similar cases in other circumstances.
\nl\nl

\bf 4. Conclusion \rm 
\nl
In Section 1 of this paper, we have introduced the concept of the completed module for a left, a right or an $A$-bimodule and in the appendix, we  see how, in a topological context, this is a completion. In Section 2, we have given a precise definition of the notion of covering and how this makes it possible to use Sweedler's notation both for coproducts and for coactions. We also have noticed that the results of the theory are not very deep, but that the examples can become quite involved. The reason for this is twofold. On the one hand, there are many different situations, all covered essentially by the same technique. On the other hand, sometimes it is needed to work successively. This is illustrated in Section 3.
\snl
As explained before, the theory is used mainly to give a meaning to expression involving Sweedler's notation. But at least as important are the cases where it is used to prove equalities of such expressions. Again, examples become sometimes complicated as before. Moreover, in this situation, an extra difficulty occurs. Often, one has to work with different coverings for the same expression and change these coverings in the course of the proof. We also have given examples where this is the case.
\snl
Finally, remark that in all cases, it is possible to avoid the use of the Sweedler notation and to give completely rigorous arguments without. However, then arguments are far less transparent. This can be seen e.g.\ in the original paper on multiplier Hopf algebra ([VD1]) where the Sweedler notation is not used at any place. 
\snl
The reader who wants to become familiar with the technique should find enough inspiration in this paper. But surely, also some personal practice will be needed to get completely acquainted with the covering technique and the use of the Sweedler notation in this context. Once this has been achieved, working with multiplier Hopf algebras is not any more difficult than working with ordinary Hopf algebras.
\nl\nl

\bf Appendix. Topological considerations \rm 
\nl
In this appendix, we will show how the completion of a module can be obtained by {\it topological methods} as the name suggest. This point of view is not important for this paper, but it certainly gives some more insight. That is why we do include it, but only in an appendix.
\snl
We again begin with the simplest case of a left $A$ module $X$ where $A$ is any algebra with a non-degenerate product. We will also assume that $A$ has right local units. We consider $X$ with the weakest topology that makes all linear functionals continuous. It is the strongest topology on $X$ making $X$ into a topological vector space. We will call this the {\it strong topology} on $X$. Now, using the module structure, we can consider a weaker topology on $X$. It is defined as the weakest topology on $X$ making all the maps $x\mapsto ax$ from $X$ to $X$ continuous for all $a\in A$ where the target space $X$ is considered with the original strong topology. Let us call this the {\it strict topology}. 
\snl
Now, consider the completed module $\overline X$, as defined in Definition 1.1 and Proposition 1.2. We can extend the strict topology to this completion by requiring that the maps $y\mapsto ay$ from $Y$ to $X$ are continuous for all $a\in A$ when again the target space $X$ is considered with the strong topology. We still call it the strict topology. 
\snl
As the topology is given by linear functionals, all these spaces are locally convex. For every locally convex space, we have the possibility to complete the space; see e.g.\ Chapter IV, Section 3 in [Y]. Then we have the following result.

\inspr{A.1} Proposition \rm
The  module $\overline X$ is the completion of the original module $X$ for the strict topology.
\snl \bf Proof\rm:
By assumption, for any finite number of elements $\{a_1, a_2, \dots, a_n\}$ in $A$, there exists an element $e\in A$ so that $a_je=a_j$ for all $j$. This provides a net $(e_\alpha)$ of elements in $A$ so that $ae_\alpha \to a$ for all $a\in A$ (in the sense that even $ae_\alpha=a$ for $\alpha$ large enough). 
\snl
Now take $x\in \overline X$ and let $x_\alpha=e_\alpha x$ for all $\alpha$. This defines a net $(x_\alpha)$ in $X$ satisfying $ax_\alpha\to ax$ for all $a$ as $ax_\alpha=ae_\alpha x$ and $ae_\alpha\to a$. Therefore, $x_\alpha\to x$ in the strict topology. This proves that $X$ is dense in $\overline X$.
\snl
We also have to show that $\overline X$ is complete. To do this, take a Cauchy net $(x_\alpha)$ in $\overline X$ with the strict topology. Take now $a\in A$ and choose $e$ so that $ae=a$. We find a Cauchy net $(ex_\alpha)$ in $X$ with the original strong topology. As $X$ is complete for this topology, we have a limit in $X$. Now define $x\in \overline X$ by $ax=a\lim ex_\alpha$. It is not hard to see that this indeed defines an element $x$ in $\overline X$ and that $x=\lim x_\alpha$. \hfill $\square$   
\einspr

The result justifies the notation and the terminology used in this paper. 
\snl
Also because of this result, it should not come as a surprise that repeating the procedure on the module $\overline X$ yields nothing new (see a remark following Notation 1.3 in Section 1).
\snl 
A similar result is of course true for the completion of a right module. Also, we can show it for an $A$-bimodule $X$. In this case, the relevant strict topology on $X$ is the weakest one so that both the maps $x\mapsto ax$ and $x\mapsto xa$ are continuous for all $a\in A$. As this topology is stronger than the two topologies considered for the left and right module structure, we see that indeed, the completion of the bimodule can be considered as sitting in the intersection of the two completions as a left, resp.\ right module (see this item in Section 1).
\nl
Next we come to the {\it topological interpretation} of {\it the extension of linear maps}. 
\snl
So as before, let $A$  have local units and let $F:X\to V$ be a linear map from the left $A$-module $X$ to a vector space $V$.
\snl
We  get the following result.

\inspr{A.2} Proposition \rm
If the variable is covered, than the map $F:X \to V$ is continuous when $X$ is considered with the strict topology and $V$ with its strongest vector space topology. The extended map (as defined in Section 2) is the unique strictly continuous extension of $F$ to the completion $\overline X$.
\snl \bf Proof\rm:
Take any $e\in A$. By assumption, the map $x\mapsto ex$ is continuous from $X$ with the strict topology to $X$ with its strong topology. Now we can compose with the linear map $F$ and then $x\mapsto F(ex)$ will be continuous from $X$ to $V$
with the strict topology on $X$. By taking $e$ so that $F=F(e\,\cdot\,)$, we find that $F$ is continuous as stated in the proposition.
\snl
As $\overline X$ is the completion of $X$ in the topological sense, we can extend the linear map $F$ to the completion by continuity. We denote this extension still by $F$. For any $x\in \overline X$, we have $x_\alpha\to x$ when $x_\alpha=e_\alpha x$ as in the proof of Proposition A.1. For $\alpha$ large enough, we get $F(e_\alpha x)=F(ex)$ and so the extension of $F$ we get here is the same as the extension obtained in Proposition 2.3. \hfill $\square$
\einspr

We have a similar result for right modules. In the case of a $A$-bimodule, the result is as follows.

\inspr{A.3} Proposition \rm Let $X$ be an $A$-bimodule and $F:X\to V$ a linear map so that the variable is covered (in the sense of Definition 2.4). Then $F$ is continuous with respect to the strict topology on the bimodule $X$ and the strong topology on the vector space $V$.
The extended map (as defined in Proposition 2.5) is the unique strictly continuous extension of $F$ to the completion $\overline X$.
\einspr

This point of view can be used when working with extended modules and extensions of maps as developed in this note.
\nl\nl

\bf References \rm
\nl
{\bf [A]} E.\ Abe: \it Hopf algebras. \rm Cambridge University Press (1977).
\snl
{\bf [De]} L.\ Delvaux: {\it The size of the intrinsic group of a multiplier Hopf algebra}. Comm.\ Alg. 31 (2003), 1499-1514. 
\snl
{\bf [De-VD-W]} L.\ Delvaux, A. Van Daele \& Shuanhong Wang: {\it Bicrossed product of multiplier Hopf algebras}.  Preprint University of Hasselt, K.U.\ Leuven and Southeast University Nanjing, {\it Version} 2 (2008).
\snl
{\bf [Dr-VD]} B.\ Drabant \& A. Van Daele: {\it Pairing and Quantum double of multiplier Hopf algebras}.  Algebras and Representation Theory 4 (2001), 109-132.
\snl
{\bf [Dr-VD-Z]} B.\ Drabant, A.\ Van Daele \& Y.\ Zhang: {\it Actions of multiplier Hopf algebras}. Commun.\ Alg.\ 27 (1999), 4117--4172.
\snl
{\bf [S]} M.\ Sweedler: {\it Hopf algebras}. Benjamin, New-York (1969). 
\snl
{\bf [VD1]} A.\ Van Daele: {\it Multiplier Hopf algebras}. Trans.\ Am.\ Math.\ Soc.\ 342 (1994), 917--932.
\snl
{\bf [VD2]} A.\ Van Daele:  {\it An algebraic framework for group duality}.  Adv.\ Math.\ 140 (1998), 323--366.
\snl
{\bf [VD-Z1]} A.\ Van Daele \& Y.\ Zhang: {\it A survey on multiplier Hopf algebras}. In 'Hopf algebras and Quantum Groups', eds. S.\ Caenepeel \& F.\ Van Oystaeyen, Dekker, New York (1998), pp. 259--309.
\snl
{\bf [VD-Z2]} A.\ Van Daele \& Y.\ Zhang: {\it Galois Theory for multiplier Hopf algebras with integrals}. Algebra and Representation Theory {\bf 2} (1999), 83-106.
\snl
{\bf [Y]} K.\ Yosida: {\it Functional analysis}. Sixth edition, Springer, Berlin (1980)
\snl

\end